\documentclass[12pt]{article}
\usepackage[english]{babel}
\usepackage{amssymb,amsmath,graphics}
\parskip\medskipamount
\newtheorem{theorem}{Theorem}[section]

\newtheorem{lemma}[theorem]{Lemma}
\newtheorem{corollary}[theorem]{Corollary}

\def\<{\langle}
\def\>{\rangle}
\newcommand{\proof}{\emph{Proof.~}}

\newcommand{\F}{\mathsf{F}}

\def\qed{{\hfill\hphantom{.}\nobreak\hfill$\Box$}}

\newcommand{\inc}{\mbox{\tt I}}
\newcommand{\ssW}{\mathsf{W}}

\begin{document}
\author{Koen Struyve\thanks{The author is supported by  the Fund for Scientific Research -- Flanders (FWO - Vlaanderen)}
}

\title{Moufang sets related to polarities in exceptional Moufang quadrangles of type $\F_4$}

\date{}
\maketitle

%\author{Koen Struyve}

%\date{}

\begin{abstract}
Departing from a Moufang set related to a polarity in an exceptional Moufang quadrangle of type $\F_4$, we construct a rank three geometry. The main property of this new geometry is that its automorphism group is identical to the one of the underlying Moufang set, providing a tool to study this Moufang set in a geometrical way. As a corollary we obtain that every automorphism of an exceptional Moufang quadrangle of type $\F_4$ stabilizing the absolute points of a polarity, also centralizes the polarity. This handles the final case of a similar result for all polarities of Moufang $n$-gons with $n \geq 3$.
\end{abstract}

\section{Introduction}

Moufang sets are rank 1 buildings satisfying the Moufang property. They are the glue in constructions of buildings with a big automorphism group (see e.g.~\cite{Ron-Tit:87})  and they play a fundamental role in the classification of certain twin buildings (see \cite{Mue:02}). The study of these groups is made harder by the lack of a suitable canonical geometry. In one of his last lectures at coll\`ege de France, Jacques Tits suggested to remedy this by defining a kind of $2$-design using the nilpotent structure of the root groups. A general problem is then to determine the collineation group of this geometry with the hope to obtain a Fundamental Theorem for the Moufang sets, i.e., the canonical situation should be that the automorphism group of the geometry is precisely the automorphism group of the Moufang set (which is that part of the automorphism group of the group that preserves the root group structure; in many cases these just coincide). This had already been done in many situations, and in the present paper we deal with the last ``Frobenius-twisted'' case, i.e., the last case where the definition of the Moufang set involves a square root of the Frobenius endomorphism. Such Moufang sets all arise from polarities in Moufang generalized $n$-gons, with $n$ even. 

The well known examples of the situation in the previous paragraph are the Suzuki groups and the Ree groups (which also exist in the finite case). The corresponding fundamental theorems are proved in \cite{hvm2} and \cite{haot}. In the present paper we deal with a much lesser known and more recent case. We consider the Moufang sets that arise from a polarity of an exceptional Moufang quadrangle of type $\F_4$. These quadrangles have been discovered by Richard Weiss in the course of the classification of Moufang quadrangles in 1997 using the root groups and commutation relations (see \cite{titsweiss}). In \cite{hvmmhl}, it is proved that these quadrangles arise from involutions in buildings of type $\F_4$ in much the same way as the exceptional quadrangles of types $\mathsf{E}_n$, $n=6,7,8$ arise from buildings of type $\mathsf{E}_n$. But a building of type $\F_4$ has a symmetric diagram, and certain ones even \emph{are} symmetric themselves. Among these, there also exist examples that allow an involutive symmetry! In some cases, these properties carry over to the exceptional Moufang quadrangles (if the involution defining these quadrangles commutes with the given symmetry), and so we obtain polarities in some exceptional Moufang quadrangles of type $\F_4$, which were only recently discovered, see \cite{hvmmhl2}. The root groups have nilpotency class 3, and so we can define canonical 2-designs in several ways. In fact, we can define a canonical geometry of rank 3 in the same way as was done for the Ree groups in \cite{haot}. In the present paper, we construct these geometries and prove the corresponding fundamental theorems.     We also list some consequences in Section~\ref{main}.

\section{Definitions and concepts}
\subsection{Moufang sets and rank one buildings}
Let $X$ be a set, with for each $x \in X$ a group $U_x$ (we name these the \emph{rootgroups}) acting on $X$ while fixing $x$. Then $(X,(U_x)_{x \in X})$ will be a \emph{Moufang set} if the following two conditions are met :
\begin{itemize}
\item For every $x \in X$, $U_x$ acts regularly on $X \backslash \{ x\}$,
\item The set of all rootgroups is normalized by the group $G^\dagger$ generated by all the rootgroups.
\end{itemize}
The group $G^\dagger$ is called the \emph{little projective group}. The \emph{full projective group} is the group of all elements of $\mathrm{Sym}(X)$ who leave the set of rootgroups invariant. 

Choose an $x \in X$, let $V_x$ be an nontrivial subgroup of $U_x$ such that $V_x$ is a normal subgroup of $G^\dagger_x$. Now we can define for each $y \in X$ an unique subgroup $V_y$ of $U_y$ as the conjugate of $V_x$ by an arbitrary element $g \in G^\dagger$ with $x^g=y$. The \emph{rank one Moufang building defined on $X$ by $(U_x)_{x \in X}$ relative to $(V_x)_{x \in X}$} is the geometry $(X,\Delta)$ where $\Delta = \{\{x\} \cup \{y^v| v \in V_x \} | x,y \in X \wedge x \neq y \}$, we will call these objects blocks. The element $x$ of a block $\{x\} \cup \{y^v| v \in V_x \}$ will be called the gnarl of that block.
\subsection{Exceptional quadrangles of type $\F_4$}
Suppose $\mathbb{K}$ a field of characteristic 2 and let $\mathbb{L}$ be a separable quadratic extension of $\mathbb{K}$. Denote by $x \mapsto \bar{x}$ the non-trivial involutory field automorphism of $\mathbb{L}$ fixing $\mathbb{K}$ pointwise. Let $\mathbb{K}'$ be a subfield of $\mathbb{K}$ containing the field $\mathbb{K}^2$ of all squares of $\mathbb{K}$ and let $\mathbb{L}'$ be the subfield generated by both $\mathbb{L}^2$ and $\mathbb{K}'$. We then have that $\mathbb{L}^2 \subseteq \mathbb{L}' \subseteq \mathbb{L}$ and because the map $x \mapsto \bar{x}$ restricted to $\mathbb{L}'$ has the fixed subfield $\mathbb{K}'$, $\mathbb{L}'$ will be a separable quadratic extension of $\mathbb{K}'$. Suppose we have two given elements $\alpha \in \mathbb{K}'$ and $\beta \in \mathbb{K}$ such that one of the following two equivalent conditions hold:
\begin{itemize}
\item $\forall u,v \in \mathbb{L}, a \in \mathbb{K}' : u \bar{u}+ \alpha v\bar{v} + \beta a= 0 \Rightarrow u=v=a=0$
\item $\forall x,y \in \mathbb{L}', b \in \mathbb{K} : x \bar{x}+\beta^2 y \bar{y}+ \alpha b^2 =0 \Rightarrow x=y=b=0$
\end{itemize}
We now identify $U_1$ and $U_3$ with the (additive) group $\mathbb{L}' \times \mathbb{L}' \times \mathbb{K}$, and $U_2$ and $U_4$ with $\mathbb{L} \times \mathbb{L} \times \mathbb{K}'$. The following relations now define the generalized quadrangle $Q(\mathbb{K},\mathbb{L},\mathbb{K}',\alpha,\beta)$ also known as the \emph{exceptional quadrangle of type $\F_4$}.
\begin{equation}
[U_1,U_2]=[U_2,U_3]=[U_3,U_4]
\end{equation}
and
\begin{align}
[(x,y,b)_1,(x',y',b')_3] =& (0,0,\alpha(x\bar{x}' + x'\bar{x} + \beta^2 (y\bar{y}' + y' \bar{y})))_2 \\
[(u,v,a)_2,(u',v',a')_4] =& (0,0,\beta^{-1} (u\bar{u}' + u'\bar{u} + \alpha (v\bar{v}' + v' \bar{v})))_2 \\
[(x,y,b)_1,(u,v,a)_4] =& (bu+\alpha(\bar{x}v +\beta y\bar{v}), bv+xu+\beta y \bar{u}, \\
& b^2 a + a \alpha (x\bar{x} + \beta^2 y\bar{y})  \nonumber  \\ 
& + \alpha (u^2 x \bar{y} + \bar{u}^2 \bar{x}y + \alpha (\bar{v}^2 xy  +v^2 \bar{x}\bar{y})))_2 \nonumber \\
& (ax+\bar{u}^2y + \alpha v^2 \bar{y},ay + \beta^{-2}(u^2 x+\alpha v^2 \bar{x}), \nonumber \\
& ab + b\beta^{-1} (u\bar{u} + \alpha v \bar{v}) \nonumber \\
& + \alpha(\beta^{-1}(xu\bar{v} + \bar{x}\bar{u}v) + y\bar{u}\bar{v} + \bar{y}uv))_3 \nonumber
\end{align}
These Moufang quadrangles were discovered by Richard Weiss in 1997 in preparation of \cite{titsweiss}. In 1999 Bernhard M\"uhlherr and Hendrik Van Maldeghem gave an algebraic interpretation to these quadrangles by constructing them out of an involution of a $\F_4$ Moufang building (\cite{hvmmhl}), furthermore they showed in \cite{hvmmhl2} that these quadrangles have polarities if and only if there is a \emph{Tits endomorphism} $\phi$ of $\mathbb{L}$ (i.e., an endomorphism such that $(x^\phi)^\phi = x^2$ for all $x$ in $\mathbb{L}$) with image $\mathbb{L}'$ such that the restriction of $\phi$ to $\mathbb{K}$ is also a Tits endomorphism with image $\mathbb{K}'$. In this case $\beta^\phi$ will be equal to $\alpha$. We let $\theta$ be the inverse of $\phi$ (defined on $\mathbb{L}'$). From now on, we will always suppose we have a polarity and we have chosen our coordinate system in such a way that this polarity takes the easy form:
\begin{align}
& (x,y,b)_1 \mapsto (\beta x^\theta, \beta y^\theta, b^{2 \theta})_4 \\
& (u,v,a)_2 \mapsto (\alpha^{-1}u^{2 \theta}, \alpha^{-1} v^{2 \theta}, a^\theta)_3 \\
& (x,y,b)_3 \mapsto (\beta x^\theta, \beta y^\theta, b^{2 \theta})_2 \\
& (u,v,a)_4 \mapsto (\alpha^{-1}u^{2 \theta}, \alpha^{-1} v^{2 \theta}, a^\theta)_1 
\end{align}

If we restrict the rootgroups $U_1,U_3$ to $\{0\} \times \{0\} \times \mathbb{K}$ and $U_2,U_4$ to $\{0\} \times \{0\} \times \mathbb{K}'$, then we obtain the Suzuki-Tits subquadrangle $\ssW(\mathbb{K},\phi)$.

\subsection{Moufang sets from exceptional quadrangles of type $\F_4$}
The absolute flags of a polarity of a Moufang quadrangle always form a Mou\-fang set. So in the case that an exceptional quadrangle of type $\F_4$ has a polarity, we also have a Moufang set named $\mathcal{M}(\mathbb{L},\mathbb{K},\alpha,\phi)$. We name the flag $\{(\infty),[\infty]\}$ of the quadrangle the element $(\infty)$ of the Moufang set, and the flag $\{((0,0,0),(0,0,0),(0,0,0)),[(0,0,0),(0,0,0),(0,0,0)]\}$ will be denoted as $[(0,0,0),(0,0,0)]$. 

The form of a generic element of the rootgroup $U_\infty$ of the element $(\infty)$ is given by (\cite{hvmmhl2}):
\begin{align}
 (x,y,a)_1 (u,v,b)_2 (\alpha^{-1}  u^{2 \theta} & +  a^{2 \theta}x + \beta^2 \bar{x}^{2 \theta} y + \alpha \beta^2 y^{2 \theta}\bar{y}, \\
 \alpha^{-1} v^{2 \theta} +  a^{2 \theta} y +  x^{2 \theta} x + \alpha y^{2 \theta} \bar{x} &,  b^\theta + a a^{2 \theta} + \alpha \beta (x^\theta \bar{x}^\theta + \alpha y^\theta \bar{y}^\theta ) +  \nonumber \\
   \alpha \beta ( \beta (y \bar{x}^\theta \bar{y}^\theta  + \bar{y}  x^\theta y^\theta ) +  &x x^\theta \bar{y}^\theta +  \bar{x} \bar{x}^\theta y^\theta) + u \bar{x}^\theta +  \nonumber \\  
 \bar{u} x^\theta + \alpha(v \bar{y}^\theta + \bar{v} y^\theta &))_3 (\beta x^\theta, \beta y^\theta, a^{2 \theta})_4 \nonumber
\end{align}
We will use a shorter notation for the above element and name it  the element $[(x,y,a)(u,v,b)]_\infty$. The image of $[(0,0,0),(0,0,0)]$ under $[(x,y,a)(u,v,b)]_\infty$ will be referred as $[(x,y,a)(u,v,b)]$. This way we have labeled all the absolute flags. This Moufang set has a sub Moufang set related to the Suzuki-Tits subquadrangle (by taking the restriction to $\{\infty\} \cup \{[(0,0,a),(0,0,b)] | a \in \mathbb{K},b \in \mathbb{K}' \}$).

\subsection{Rank 1 Moufang buildings constructed out of  \\ $\mathcal{M}(\mathbb{L},\mathbb{K},\alpha,\phi)$}
As $U_\infty$ has nilpotency class 3, we have 2 obvious candidates for choosing $V_\infty$, namely $[U_\infty,U_\infty]=\{[(0,0,0),(k,l,m)]_\infty | k,l \in \mathbb{L}, m \in \mathbb{K}' \}$ (subsets constructed this way we name \emph{spheres}) and $[U_\infty,[U_\infty,U_\infty]]=\{[(0,0,0), $ $(0,0,m)]_\infty |  m \in \mathbb{K}' \}$ (we name these \emph{circles}). We obtain a rank 3 geometry $\Omega = (X,Y,Z,\inc)$ of the points in the Moufang set ($X$), the spheres ($Z$), circles ($Y$) and as incidence relation $\inc$ containment. The circles who are completely in the restriction to $\{\infty\} \cup \{[(0,0,a),(0,0,b)] | a \in \mathbb{K},b \in \mathbb{K}' \}$ form the blocks of a similar geometry constructed from the Moufang set related to the Suzuki-Tits quadrangle, of which a similar study has been done by Hendrik Van Maldeghem in \cite{hvm2}. We will reuse several techniques used in the latter paper here.

\section{Statement of main results}\label{main}
In the previous section we have constructed the geometry $\Omega$, our question is now what is the automorphism group of this geometry? We obtained a stronger result, in the form that both the automorphism group of $\Omega$ and the truncated geometry  $(X,Z,\inc)$, which consists of only the points and the spheres, are canonically isomorphic to the subgroup of the automorphism group of the exceptional quadrangle of type $\F_4$ centralizing the polarity and the full projective group of the Moufang set $\mathcal{M}(\mathbb{L},\mathbb{K},\alpha,\phi)$.

%and the truncated geometries $(X,Y,\inc)$ and $(X,Z,\inc)$. In the case of $\Omega$ and $(X,Z,\inc)$ we have determined the automorphism group, in the case of $(X,Y,\inc)$ we only have a partial result.
\begin{theorem} \label{theorem:t1}
The automorphism group of $(X,Z,\inc)$ is isomorphic with the subgroup of the automorphism group of the exceptional quadrangle of type $\F_4$ centralizing the polarity and the full projective group of the Moufang set $\mathcal{M}(\mathbb{L},\mathbb{K},\alpha,\phi)$.
\end{theorem}
\begin{corollary} \label{collorary:c1}
The automorphism group of $\Omega$ is isomorphic with the subgroup of the automorphism group of the exceptional quadrangle of type $\F_4$ centralizing the polarity and the full projective group of the Moufang set $\mathcal{M}(\mathbb{L},\mathbb{K},\alpha,\phi)$.
\end{corollary}
This knowledge about the automorphisms centralizing the polarity has the following interesting corollary:
\begin{corollary}
Each automorphism of an exceptional quadrangle of type $\F_4$ stabilizing the set of absolute points of a certain polarity of the quadrangle, also stabilizes the set of absolute lines and will centralize that polarity.
\end{corollary}
As the other cases have been handled before (the case of the Suzuki-Tits quadrangle in~\cite{hvm1}, and the Ree hexagon in~\cite{haot}), a more general result is now immediate.
\begin{corollary}
Each automorphism of a Moufang $n$-gon with a polarity stabilizing the set of absolute points of that polarity, also stabilizes the set of absolute lines and centralizes that polarity, except if either $n=3$, the projective plane is Pappian, the characteristic of the underlying field is $2$, and the polarity is not Hermitian (i.e., there is no twisting field automorphism); or if $n=4$ and the generalized quadrangle is the smallest symplectic quadrangle $\ssW(2)$.
\end{corollary}

The polarities of the exceptions when $n=3$ are sometimes called `pseudo-polarities'. The set of absolute points of such a pseudo-polarity is the set of points on a line $L$, and one additionally needs to fix the unique point (not on $L$) through which every line is an absolute line (and this point can be chosen arbitrarily). 

\section{Proof of main results}
\subsection{Automorphism group of $(X,Z,\inc)$}
We denote this geometry with $\Omega'$ to make the notations easier. It is easy to see that each automorphism of the quadrangle centralizing the polarity, gives canonically an element of the full projective group, which on its turn gives an automorphism of $\Omega'$. The aim of the proof is to show that we can return from an automorphism of $\Omega'$ to an automorphism of the quadrangle centralizing the polarity, in such a way that applying these 3 consecutive group morphisms brings you back to the original element. 

The sphere with gnarl $(\infty)$ containing $[(0,0,0),(0,0,0)]$ is given by the union of the element $\{ (\infty) \}$ with the orbit of $[U_\infty,U_\infty]$ on $[(0,0,0),(0,0,0)]$, which will be $\{[(0,0,0),(k,l,m)] | k,l \in \mathbb{L}, m \in \mathbb{K}' \}$. The first lemma gives us a geometric interpretation of the spheres.

\begin{lemma} \label{lemma:constr}
The collection of absolute flags in the quadrangle forming the sphere with gnarl $\{p,L\}$ and containing $\{q,M\}$ are those absolute flags of which the point is collinear with the projection of $q$ onto $L$.
\end{lemma}
\proof
If we look at the specific case of the flags $\{p,L\}=(\infty)$ and  $\{q,M\}=[(0,0,0),(0,0,0)]$, we can see this is true by using coordinates within the quadrangle. Because the little projective group of the Moufang set, which also acts on $\Omega$, is constructed as a subgroup of the automorphism group of the quadrangle which is 2-transitive on the absolute flags and preserves collinearity this property will be true for all choices of points in $X$.
\qed

This lemma allows us to calculate several spheres, the ones we need are listed in the first appendix. The \emph{derived geometry} $\Omega'_{(\infty)}$ is the geometry formed by the points in $X\backslash{(\infty)}$  and the intersections of the spheres through $(\infty)$ with $X\backslash{(\infty)}$ (called lines). The lines of $\Omega'_{(\infty)}$ which have $(\infty)$ as gnarl will be called the vertical lines, the others the non-vertical lines. The group $U_\infty$ acts transitively on both the set of the vertical lines and the non-vertical ones ($U_\infty$ does not map vertical lines to non-vertical lines or vice versa). The main question is if the stabilizer of $(\infty)$ in $\Omega$ has one or two orbits on all the spheres through $(\infty)$. We will start under the assumption that there is only one orbit (this will be denoted by \cal{[T]}) and try to examine the structure of $\Omega'_{(\infty)}$ in more detail.

We first need an additional definition. A \emph{net} is a geometry $(P,B,\inc)$ consistings of points ($B$), lines ($B$) and an incidence relation ($\inc$) such that the following 2 conditions are met:
\begin{itemize}
\item
Two lines have at most one point in common.
\item
For each line $L \in B$ and point $x \in P$ not incident with $L$, there is exactly one line incident with $x$ not intersecting $L$.
\end{itemize}
These properties imply that $B$ is partitioned in parallel classes, which are sets of mutually non-intersecting lines such that each point is incident with exactly one line of the parallel class.

\begin{lemma} 
If condition \cal{[T]} holds then $\Omega'_{(\infty)}$ will be a net.
\end{lemma}
\proof
By condition \cal{[T]} one can always choose one line to be a vertical line when checking these properties.  We first compute in how many points 2 such lines will intersect :
\begin{itemize}
\item Two vertical lines will never intersect as they form differents orbits of $[U_\infty,U_\infty]$,
\item One vertical line and one non-vertical line : by using transitivity of $U_{\infty}$ on the non-vertical lines, we can suppose that the non-vertical line is $B_{[(0,0,0),(0,0,0)]}$. The intersection point with a generic vertical block $B_{(x,y,a)}$ is $[(x,y,a),(0,0,0)]$, so these always intersect in exactly one point.
\end{itemize}
This also proves that the vertical lines form a parallel class for which the second axiom of nets holds, by transitivity we obtain that $\Omega'_{(\infty)}$ is a net. 
\qed

We remark that under assumption of $[T]$, the parallel class of $B_{[(k,l,m),(u,v,b)]}$ is given by all blocks of the form $B_{[(k,l,m),(u',v',b')]}$ (this is easily verified as the exact form in the appendix implies that they do not intersect if $(u,v,b) \neq (u',v',b')$).

The non-identity elements of $[U_\infty,U_\infty]$ will fix $(\infty)$ and all the vertical lines of $\Omega'_{(\infty)}$ and acts freely on the points on these lines. If \cal{[T]} holds then there would be similar automorphisms of $\Omega$ for other parallel classes than the vertical lines. The following lemma excludes this possibility. 
\begin{lemma}
If $\Omega'_{(\infty)}$ is a net then no automorphism of $\Omega$ fixes the point $(\infty)$ and all lines of $\Omega'_{(\infty)}$ parallel with $B_{[(0,0,0),(0,0,0)]}$, acts freely on the points of these lines and maps $[(0,0,0),(0,0,0)]$ to $[(0,0,1),(0,0,0)]$.
\end{lemma}
\proof
Suppose it does exists and let $\tau$ be such an automorphism.

If a line $L$ of $\Omega'_{(\infty)}$ would be mapped to a line $L^\tau$ of another parallel class then the intersection point $L \cap L^\tau$ would be fixed because all the lines parallel to $B_{[(0,0,0),(0,0,0)]}$ are fixed. But we want $\tau$ to act freely on the points of these lines so we have that all parallel classes are stabilized.

Because parallel classes are stabilized we have that $B_{(0,0,0)}^\tau = B_{(0,0,1)}$. As a consequence of this we have that the point $[(0,0,0),(u,v,b)]$ will be mapped to $[(0,0,1),(u,v,b)]$ (as both points have to be on the same block $B_{[(0,0,0),(u,v,b)]}$ parallel to $B_{[(0,0,0),(0,0,0)]}$). This implies that the block $B_{[(0,0,1),(u,v,b+1)]}$ through the point $[(0,0,0),(u,v,b)]$ will be mapped to the block $B_{[(0,0,1),(u,v,b)]}$. Considering the intersections of these 2 lines with the fixed line $B_{[(0,0,0),(0,0,0)]}$, we have that $[(\alpha^{-1} u^{2 \theta},\alpha^{-1} v^{2 \theta},b^\theta),(0,0,0)]$ is mapped to $[(\alpha^{-1} u^{2 \theta},\alpha^{-1} v^{2 \theta},b^\theta+1),(0,0,0)]$, or if we rename the variables: $[(x,y,a),(0,0,0)]^\tau=[(x,y,a+1),(0,0,0)]$. The generic point $[(x,y,a),(u,v,b)]$ lies on the lines $B_{[(0,0,0),(u,v,b)]}$ and $B_{(x,y,a)}$, these lines are mapped to $B_{[(0,0,0),(u,v,b)]}$ and $B_{(x,y,a+1)}$, so we have that $[(x,y,a),(u,v,b)]^\tau=[(x,y,a+1),(u,v,b)]$.

If we first apply the automorphism $[(0,0,1),(0,0,0)]_{(\infty)}$ followed by $\tau$, we obtain an automorphism $\tau'$ which maps $[(x,y,a),(u,v,b)]$ to $[(x,y,a),(u+\beta x^\theta,v+x^\theta,b+a^{2\theta})]$. 

All that is left to show is that the image under $\tau'$  of the sphere $S$ with gnarl $[(0,0,0),(0,0,0)]$ through $[(0,0,1),(0,0,0)]$ isn't a sphere, giving us a contradiction. The sphere $S$ will contain the circle $C$ with gnarl $[(0,0,0),(0,0,0)]$ through $[(0,0,1),(0,0,0)]$, which on its turn contains the point $[(0,0,1),$ $(0,0,1)]$ (see the second appendix). Hence $S^{\tau'}$ will contain the points $[(0,0,0),$ $(0,0,0)]$, $[(0,0,1),(0,0,0)]$ and $[(0,0,0),(0,0,1)]$. As $\Omega'_{(\infty)}$ is a net two different spheres of $\Omega'$ will have at most two points in common. This implies $S^{\tau'}$ will be the sphere through $[(0,0,1),(0,0,0)]$ with gnarl $[(0,0,0),(0,0,1)]$ (because these three points points lay on the circle containing $[(0,0,1),(0,0,0)]$ with gnarl $[(0,0,0),(0,0,1)]$, which is contained in the sphere). Moreover, $\tau'$ stabilizes the set $\{[(0,0,a),(0,0,b)] | a \in \mathbb{K},b \in \mathbb{K}' \}$ so that the circle $C$ will be mapped to the circle with gnarl $[(0,0,1),(0,0,0)]$ through $[(0,0,0),(0,0,1)]$. That this is impossible is shown by Hendrik Van Maldeghem in lemma 4 of \cite{hvm2}.
\qed

All this has proved that the condition \cal{[T]} is false, implying that the gnarls of the spheres are uniquely defined, and that the gnarl of a given sphere can be recognized solely by looking at the properties of the geometry $\Omega'$. The last step of our proof is to show that we can reconstruct the exceptional Moufang quadrangle of type $\F_4$ out of $\Omega'$ and see that the automorphisms of $\Omega'$ lift to the desired automorphisms of the quadrangle. The construction of the quadrangle goes as follows :
\begin{itemize}
\item Points : 
\begin{itemize}
\item points $x$ of $\Omega'$, notation as point of the quadrangle : $x_p$
\item spheres $B$ of $\Omega'$, notation as point of the quadrangle : $B_p$
\end{itemize}
\item Lines : 
\begin{itemize}
\item points $x$ of $\Omega'$, notation as line of the quadrangle : $x_l$
\item spheres $B$ of $\Omega'$, notation as line of the quadrangle : $B_l$ 
\end{itemize}
\item Incidence :
\begin{itemize}
\item $x_p \; I \; y_l \Leftrightarrow x=y$
\item $x_p \; I \; A_l \Leftrightarrow A_p \; I \; x_l \Leftrightarrow $ the gnarl of $A$ is $x$
\item $A_p  \; I \; B_l \Leftrightarrow $ the gnarl of $A$ is contained in $B$, the gnarl of $B$ is contained in $A$ and these 2 gnarls are different
\end{itemize}
\end{itemize}
That this construction gives back the exceptional Moufang quadrangle of type $\F_4$ follows from lemma \ref{lemma:constr}. The polarity of the quadrangle in this construction comes down to interchanging the subscripts $_p$ and $_l$. It is also easily seen that an automorphism of $\Omega'$ lifts naturally to an automorphism of the exceptional quadrangle of type $\F_4$ centralising the polarity in the same way that an automorphism of the exceptional quadrangle of type $\F_4$ centralising the polarity implies an automorphism of $\Omega'$. This proves theorem \ref{theorem:t1}.
\subsection{Proof of the corollaries}

Corollary \ref{collorary:c1} follows from theorem \ref{theorem:t1} because an automorphism of the exceptional quadrangle of type $\F_4$ centralising the polarity will also act well-defined on the sets of the absolute flags forming the circles and thus implies an automorphism of $\Omega$.

The second corollary follows from the fact that due to lemma \ref{lemma:constr} an automorphism $\eta$ of the quadrangle which stabilizes the absolute points (well-defined action on $X$) and preserves collinearity (the lemma tells us that if this is preserved blocks will be mapped to blocks) also implies an automorphism of $(X,Z,\inc)$. But we have proven that such an automorphism of $(X,Z,\inc)$ comes from an automorphism of the quadrangle centralising the polarity. So $\eta$ stabilizes the set of absolute lines and will centralize the polarity, which proves the second corollary.

%\subsection{Automorphism group of $(X,Y,\inc)$}
%Here we have only a weaker result saying the gnarls of the circles are preserved by automorphisms of $(X,Y,\inc)$. If we look to the formulas for circles through $(\infty)$ in the second appendix, we see that all the points different from $(\infty)$ of one particular circle all have a constant first, second, fourth and fifth coordinate. Moreover, this way we can see that the derived geometry of $(X,Y,\inc)$ in a point will be disconnected, and so that each of the connected components are the points of a  rank 1 Suzuki-Tits sub Moufang building. This way we can interprete the circles as blocks of that subbuilding and use the result that the blocks of a rank 1 Suzuki-Tits Moufang building have unique blocks (proved by Hendrik Van Maldeghem in \cite{hvm2}) to prove that the circles have unique gnarls. This proves the second main result.

\appendix

\section{Coordinates of points of certain spheres}
Spheres going through $(\infty)$ :
\begin{itemize}
\item spheres with gnarl $(\infty)$ going through the point $[(x,y,a),(u,v,b)]$ :
\begin{equation}
B_{(x,y,a)} = \{(\infty)\} \cup \{[(x,y,a),(k,l,m)] | k,l \in \mathbb{L}, m \in \mathbb{K}' \}
\end{equation} 
\item spheres with gnarl $[(x,y,a),(u,v,b)]$ (denoted by $B_{[(x,y,a),(u,v,b)]}$) :
\begin{align}
B_{[(0,0,0),(u,v,b)]} =  &\{(\infty)\} \cup \{[(k,l,m),(u,v,b)] | k,l \in \mathbb{L}', m \in \mathbb{K} \}  \\
B_{[(0,0,1),(u,v,b)]} =  &\{(\infty)\} \cup \{[(k,l,m+1),(u+\beta k^\theta ,v+ \\
 &\beta l^\theta,b+m^{2\theta})] | k,l \in \mathbb{L}', m \in \mathbb{K} \}
\end{align}

\end{itemize}
\section{Coordinates of points of certain circles}
Circles through $(\infty)$ :
\begin{itemize}

\item circles with gnarl $(\infty)$ going through the point $[(x,y,a),(u,v,b)]$ :
\begin{equation}
\{(\infty)\} \cup \{[(x,y,a),(u,v,k)] | k \in \mathbb{K}' \}
\end{equation} 
\item circles with gnarl $[(x,y,a),(u,v,b)]$ :
\begin{equation}
\{(\infty)\} \cup \{[(x,y,a+k),(u,v,b+a^2 k^{2\theta} + k^{2\theta}\alpha (x\bar{x} +\beta^2 y\bar{y}) )] | k \in \mathbb{K}' \}
\end{equation} 
\end{itemize}
The circle with gnarl $[(0,0,0),(0,0,0)]$ through $[(0,0,1),(0,0,0)]$ :
\begin{equation}
\{[(0,0,0),(0,0,0)]\} \cup \{[(0,0,\frac{x}{1+x+x^{2\theta}}),(0,0,\frac{1}{1+x^{2\theta}+x^2})] | x\in \mathbb{K} \}
\end{equation}
The circle through $[(0,0,1),(0,0,0)]$ with gnarl $[(0,0,0),(0,0,1)]$ :
\begin{equation}
\{[(0,0,0),(0,0,0)]\} \cup \{[(0,0,\frac{1}{1+x+x^{2+2\theta}}),(0,0,\frac{x^{2 \theta}}{1+x^{2\theta}+x^{2+4\theta}})] | x\in \mathbb{K} \}
\end{equation}


\begin{thebibliography}{9}
\bibitem{haot}
Fabienne Haot, Koen Struyve and Hendrik Van Maldeghem. Ree Geometries, \emph{Forum Math.}, accepted
\bibitem{Mue:02} Bernhard M\"uhlherr, `Locally split and locally finite twin buildings of 2-spherical type', \emph{J. Reine Angew. Math.} \textbf{511} (1999), 119--143.
\bibitem{hvmmhl}Bernhard M\"uhlherr and Hendrik Van Maldeghem, Exceptional Moufang quadrangles of type $\mathrm{F}_4$, \emph{Canad. J. Math.} \textbf{51} (1999), 347--371
\bibitem{hvmmhl2}
Bernhard M\"uhlherr and Hendrik Van Maldeghem, Moufang Sets from Groups of Mixed Type, \emph{J. Algebra} \textbf{300}  (2006),  no. 2, 820--833.
\bibitem{Ron-Tit:87} 
Mark Ronan and Jacques Tits, Building buildings, {\em Math.\ Ann.}\ {\bf 278}, 291--306.
\bibitem{titsweiss}
Jacques Tits and Richard Weiss, \emph{Moufang Polygons}, Springer Monographs in Mathmematics, Springer-Verlag, 2002
\bibitem{hvm1}
Hendrik Van Maldeghem, \emph{Generalized Polygons}, Birkh\"auser Verlag, 1998
\bibitem{hvm2}
Hendrik Van Maldeghem, Moufang lines defined by (generalized) Suzuki groups, \emph{European J. Combin.} \textbf{28} (2007), no. 7, 1878--1889. 
\end{thebibliography}
\end{document}